\documentclass[12pt]{amsart}
\usepackage{amssymb}
\usepackage{latexsym}
\usepackage{amsfonts}
\usepackage{amsmath}
\usepackage{amscd}
\usepackage{quoting}

\newtheorem{theorem}{Theorem}[section]

\theoremstyle{plain}

\newcommand{\sy}[1]{{\sf S}_{#1}}

\newcommand{\soc}{\operatorname{Soc}}

\newcommand{\twr}[3]{#1\, {\sf twr}_{#2}\,#3}
\renewcommand{\wr}{\,{\sf wr}\,}

\newcommand{\aut}[1]{{\sf Aut}\,{#1}}
\newcommand{\out}[1]{{\sf Out}(#1)}
\newcommand{\inn}[1]{{\sf Inn}\,#1}

\newcommand{\norm}[2]{{\mathbb N}_{#1}\left(#2\right)}

\renewcommand{\leq}{\leqslant}
\renewcommand{\geq}{\geqslant}

\font\tenvr=cmmi10 scaled 1600
\textfont"F=\tenvr

\title[L.\ G.\ Kov\'acs on permutation groups]{The contribution of L.\ G.\ Kov\'acs to the theory of permutation groups}
\author{Cheryl E. Praeger}
\address{School of Mathematics and Statistics\\
The University of Western Australia\\
35 Stirling Highway 6009 Crawley\\
Western Australia\\
 cheryl.praeger@uwa.edu.au\\
www.maths.uwa.edu.au/$\sim$praeger}
\author{Csaba Schneider}
\address{Departamento de Matem\'atica\\
Instituto de Ci\^encias Exatas\\
Universidade Federal de Minas Gerais\\
Av.\ Ant\^onio Carlos 6627\\
Belo Horizonte, MG, Brazil\\
csaba@mat.ufmg.br\\
 www.mat.ufmg.br/$\sim$csaba}

 \thanks{The second author was supported by the
 research projects 302660/2013-5 (CNPq, Produtividade em Pesquisa),
 475399/2013-7 (CNPq, Universal), and 
 the APQ-00452-13 (Fapemig, Universal).
He also enjoyed the hospitality of the School of Mathematics and Statistics of the University of Western Australia, while a part of this paper was 
being written.}
\begin{document}
%\maketitle
\begin{abstract}
The work of L.\ G.\ (Laci) Kov\'acs (1936--2013) 
gave us a deeper understanding of permutation groups, especially
in the O'Nan--Scott theory of primitive groups. We review his 
contribution to this field.
\end{abstract}

\dedicatory{Dedicated to the memory of our friend Laci Kov\'acs}

\date{\today}
\maketitle

\section{The context of Laci's work on permutation groups}\label{sect1}
The main contributions made by L.\ G.\ Kov\'acs (`Laci' to his friends and 
colleagues) in the theory of permutation groups revolve around 
groups acting primitively  and their point stabilizers, which are 
maximal subgroups.

The classification of finite simple groups was, undoubtedly, one of 
the greatest research projects ever undertaken in mathematics.
It changed  the face of finite group theory, and several important results 
on permutation groups grew from the classification.
Laci's work contributed significantly to one such result, 
namely the O'Nan--Scott Theorem for primitive permutation groups.
%% He, independently with Aschbacher and Scott, corrected the 
%% mistakes in its early versions, made its statement more precise,
%% and showed how to use it 
%% to obtain information about the structure of finite groups, for
%% instance, to describe maximal subgroups.

The origins of the  O'Nan--Scott Theorem can be traced back to 
Jordan's description of the maximal subgroups of the symmetric groups in his
{\em Trait\'e}~\cite{jordan}.  The modern form of the theorem dates
back  to 1978.  By then  it  was clear  that the  announcement of  the
classification  of the  finite simple  groups was  imminent,  and many
mathematicians  began  to   consider  how  this  classification  might
influence work  in other parts  of mathematics, especially  in finite
group theory.  With different applications in mind,  Leonard Scott and
Michael  O'Nan  independently  proved     versions  of  this  theorem,
pin-pointing  various r\^oles for  the simple  groups in  describing the
possible kinds of finite  primitive permutation groups. Both Scott and
O'Nan  brought  papers  containing  their  result to  the  Santa  Cruz
Conference of  Finite Groups in  1979~\cite{santacruz}, 
in which Laci also participated.
%According to Scott~\cite{scottpage}, O'Nan commented, after
%seeing  Scott's   manuscript:    ``Same  damn
%theorem!''. 
As was kindly pointed out to us by the referee, 
a preliminary version of the conference proceedings included
the papers by  both Scott and O'Nan, while 
the final version only contained Scott's paper,
which stated the theorem
in an  appendix~\cite{sco:rep}.
%O'Nan's paper  was not published  in the 
%proceedings; 
%only Scott's  paper was published in the proceedings 
%(acknowledging that
%O'Nan  had  obtained the  theorem  independently). 
%Several  colleagues
%brought  copies of  O'Nan's manuscript  to  Australia to  give to  the
%first author, and  she was disappointed that it  was not published in
%the proceedings.

These early versions of the O'Nan--Scott Theorem were presented as 
characterizations of the maximal subgroups of the finite alternating
and symmetric groups, 
as well as structural descriptions for finite primitive permutation groups
 (see the theorems stated on pages 329 and 328 of~\cite{sco:rep}, respectively).
They were unfortunately not adequate for the latter purpose since 
the statement of the Theorem on page 328 of~\cite{sco:rep} 
failed to identify primitive groups of twisted wreath type 
(that is, those with a unique  minimal normal subgroup which is non-abelian
and regular).
In addition,  the theorems on pages 328--329 in~\cite{sco:rep} 
erroneously claimed that that 
the number of simple factors of the unique minimal normal subgroup
of a primitive group of simple diagonal type
had to be a prime number. 

These oversights were quickly recognized.
They were first rectified by Aschbacher and Scott~\cite[Appendix]{aschsc} 
and, independently of their work, 
by Laci~\cite{kovacs:max}. In particular, Laci identified the `missing type' of primitive groups: these are groups which, in his later work,
were characterized as  twisted wreath products.
%The first
%people to  notice these mistakes   were Michael  Aschbacher
%and Laci. Aschbacher with
%Scott  published  a corrected  and  expanded  version  of the  
%theorem~\cite[Appendix]{aschsc}. Independently of their work, \cite{kovacs:max}
%also identified the missing type of primitive groups: these are groups which in% his later work
%were characterized as  twisted wreath products.
%It was noticed after Corollary~3.3 of~\cite{grosskov}
%that a primitive group with a regular non-abelian minimal normal subgroup
% has the  structure of a  twisted wreath product.
%This observation
%has lent its name to the groups of
%this  type.  
A  short, self-contained  proof of  the O'Nan--Scott Theorem  
was then
published  by  Liebeck,  Saxl, and the first author~\cite{lps:onan}.  
%in  particular  proving
%additional properties of the twisted wreath types~\cite{lps:onan}.

The O'Nan--Scott Theorem has several interpretations depending on the 
context in which it is used. The theorem is most commonly interpreted, 
for instance in its applications in permutation group theory and 
algebraic combinatorics, as
a structure theorem describing finite primitive permutation groups.
It classifies such groups into several classes each of which
is characterized by the structure of the socle viewed as a permutation group. 
%The theorem is most 
%commonly used  in this form in its applications in permutation group 
%theory and algebraic combinatorics.
The fact that
maximal subgroups in finite groups are point stabilizers for primitive
actions leads to a second interpretation of 
the O'Nan--Scott Theorem:
 it can be considered as a characterization of (core-free) 
maximal subgroups of finite groups. 
Much of Laci's work related to the O'Nan--Scott Theorem was concerned with
giving a detailed description of maximal subgroups. His 
methodology relies, in large part, on the concept of induced extensions 
introduced in his 
%
%The  first major paper  of Laci that is related to the O'Nan--Scott Theorem,
%though without explicit reference,
%is the  
1984 paper~\cite{grosskov} jointly  written with
Fletcher  Gross. 
%% They  consider finite  groups  $G$ possessing  a
%% normal subgroup $M$ which is a direct product of $G$-conjugate subgroups.
%% They show that every such group can be constructed from a smaller group $A$
%% using a technique that they call {\em induced extension}.
%% They also describe 
%% an explicit correspondence between the subgroup lattice of $G$ and the 
%% subgroup lattice of $A$.
%% We discuss these results  in more detail
%% in Section 2. 
Based on this concept, he gave in~\cite{kovacs:max} a general description of maximal subgroups
of finite groups $G$ having a normal subgroup that is a direct product
of subgroups forming a conjugacy class of $G$. 

In his papers on permutation groups, Laci was often concerned with 
counting a class of mathematical objects, for instance maximal subgroups 
(as in~\cite{kovacs:max}) or wreath decompositions (in~\cite{kov:decomp}) 
of a group. His underlying philosophy is beautifully explained in the 
introduction of~\cite{kov:decomp}.
\begin{quotation}
{\small Counting the number of mathematical objects of a certain kind is often 
undertaken as a test problem (``if you can't count them, you don't really 
know them''): not so much because we want the answer, but because the 
attempt focuses 
attention on gaps in our understanding, and the eventual proof may 
embody insights beyond those which are capable of concise expression 
in displayed theorems.\quad \cite[p.~255]{kov:decomp}}
\end{quotation}

We describe Laci's  most important results in this field using
modern terminology for the O'Nan--Scott Theorem, namely we use a version
of the theorem that classifies finite primitive permutation groups
into eight types: Holomorph of Abelian (HA),
Almost Simple (AS), Simple Diagonal (SD), Compound Diagonal (CD), 
Holomorph of a Simple group (HS), Holomorph of a Compound group (HC), 
Product Action (PA), and Twisted Wreath (TW) --- with two-letter abbreviations as indicated.
This `type subdivision' was originally suggested by Laci and was first used 
in  the first author's paper~\cite{prae:incl} where a detailed description of each type is given. 
The paper~\cite{prae:incl} analyses
the possible inclusions between finite primitive groups in terms of these types, and was written in
close consultation with Laci who at the time was working on the related
article~\cite{kov:blowup}.

While Sections 2--4 focus on Laci's contributions around the O'Nan--Scott Theorem, 
in the final section we mention several other papers by Laci on permutation groups, including his papers 
on minimal degree of a group, and on abelian quotients of permutation groups.

%% In a second major paper, the induced extension is used to characterize
%% maximal subgroups of groups $G$ that satisfy the conditions
%% of the previous paragraph. He focuses on the case where 
%% $M$ is the direct product of finitely many simple groups permuted
%% by $G$.  The paper gives  a
%% succinct description  of the possible structures  of point stabilizers
%% of faithful  primitive permutation representations of a  given group $G$
%% with a non-simple, non-abelian  minimal normal subgroup, and counts the
%% number of equivalence classes of such representations of $G$. We discuss
%% this in Section~3.

\section{Induced extensions}

A central theme of Laci's papers in the 1980's was the study of 
groups that contain a normal subgroup $M$ which is an unrestricted 
direct product
$M=\prod_i M_i$ in such a way that the set of the $M_i$ is closed 
under conjugation by $G$. It was his typical technique to consider, as  building blocks in such
direct products, the normal subgroups $\prod_{j\neq i}M_j$ of $M$ 
instead of the direct factors $M_i$, and we will use this view here.
Laci usually worked under the following assumption:

\begin{quotation}
{\small  (*) $G$ is a group, $M$ is a normal subgroup of $G$, $\{K_i\}$ is a 
$G$-conjugacy class containing normal subgroups of $M$ in such a way
that the natural homomorphism from $M$ into $\prod_i M/K_i$ is an 
isomorphism. Let $K$ be a fixed element of the set $\{K_i\}$ and 
let $N=\norm GK$.}
\end{quotation}

Sometimes hypothesis (*) was relaxed to require only that the set $\{K_i\}$
was a $G$-invariant set of normal subgroups of $M$ (that is, a union
of $G$-conjugacy classes).
Groups that satisfy (*) occur naturally among finite
primitive permutation groups and Laci's results have indeed found several
deep applications in permutation group theory. 
Many of the classes of the O'Nan--Scott Theorem 
are characterized by the fact that
their groups contain a non-abelian minimal normal subgroup $M=T^k$ which is the 
direct product of finite simple groups all isomorphic to a group $T$. 
Such groups clearly satisfy (*).
%Further
%$M$ admits a direct decomposition $M=L^\ell$ in such a way that a 
%point stabilizer $M_\omega$ can also be written as the corresponding
%direct product $M_\omega=(M_\omega\cap L)^\ell$. We may assume, in many cases t%he 
%direct decompositions $T^k$ and $L^\ell$ coincide, however, there is a case, 
%when the stabilizer $M_\omega$ is a direct product of several diagonal 
%subgroups in $M$, and in this case  $T^k$ is a proper refinement of $L^\ell$. 

% 
The main tool in the investigation of the groups in 
which~(*) holds is the induced extension, introduced 
in~\cite[Section~3]{grosskov}. The construction takes as input a group homomorphism
$\alpha:A\rightarrow B$, and outputs a subgroup $G$ of the wreath product
$W=A\wr P$, where $P$ is the permutation group induced by the multiplication action of $B$ on the right coset space $I:=[B:A\alpha]$.
The group $G$ is called the {\em induced extension defined by $\alpha$}.  
To explain the construction, we use the following convention (often used by 
 Laci himself). For sets 
$X$ and $Y$, we denote by $X^Y$ the set of functions 
$Y\rightarrow X$. If $X$ is a group, then $X^Y$ is also a group
under pointwise multiplication, and is isomorphic to 
the (unrestricted) direct product of $|Y|$ copies of $X$.

% The main tool in the investigation of the groups in 
% which~(*) holds
% is the {\em induced extension} defined in~\cite[Section~3]{grosskov}. 
% To explain the construction, we use a convention often used by 
% Laci himself and regard direct products as functions. More precisely, if 
% $X$ and $Y$ are sets, we denote by $X^Y$ the set of functions 
% $Y\rightarrow X$. If $X$ is a group, then $X^Y$ can be turned into a group
% using pointwise multiplication and it will be isomorphic to 
% the (unrestricted) direct product of $|Y|$ copies of $X$.

% Let $A$ and $B$ be groups and 
% assume that $\alpha:A\rightarrow B$ is a homomorphism.
% Let $C=A\alpha$, and let $P$ be the permutation group
% induced by $B$ on the set $I$ of right cosets  
% of $C$ by multiplication. It
% is well-known that there exists an embedding $\lambda$ of 
% $B$ into the wreath product
% $C\wr P=C^I\rtimes P$ (see the discussion on~\cite[p.~136]{grosskov}).
% Further, there is an epimorphism 
% $\bar\alpha:A\wr P\rightarrow C\wr P$ whose kernel is $M=(\ker\alpha)^I$. 
% %The configuration is illustrated by the following diagram:
% %\begin{equation}\label{diag1}
% %\begin{CD}
% %@.B\\
% %@.@VV\lambda V\\
% %A\wr P @>\bar\alpha>> C\wr P.
% %\end{CD}
% %\end{equation}
% As $\lambda$ is a monomorphism,  there is a unique ``smallest'' group $G$
% such that suitable homomorphisms make the following pull-back 
% diagram commutative:

Let $\alpha, I, P$ be as in the previous paragraph, and set $C=A\alpha$.
It is well-known that there exists an embedding $\lambda$ of 
$B$ into the wreath product
$C\wr P=C^I\rtimes P$ (see the discussion on~\cite[p.~136]{grosskov}).
Further, there is an epimorphism 
$\bar\alpha:A\wr P\rightarrow C\wr P$ whose kernel is $M=(\ker\alpha)^I$. 
%The configuration is illustrated by the following diagram:
%\begin{equation}\label{diag1}
%\begin{CD}
%@.B\\
%@.@VV\lambda V\\
%A\wr P @>\bar\alpha>> C\wr P.
%\end{CD}
%\end{equation}
As $\lambda$ is a monomorphism,  there is a unique `smallest' group $G$
such that suitable homomorphisms make the following pull-back 
diagram commutative:
\begin{equation}\label{diag2}
\begin{CD}
G@>>>B\\
@VVV @VV\lambda V\\
A\wr P @>\bar\alpha>> C\wr P.
\end{CD}
\end{equation}
Further, as $M=(\ker\alpha)^I$ is the kernel of $\bar\alpha$, we obtain the short
exact sequence
$$
\begin{CD}
1 @>>> M @>>> A\wr P @>\bar\alpha >> C\wr P @>>> 1.
\end{CD}
$$
Combining the pull-back diagram with this short exact sequence we obtain
the pull-back exact sequence
$$
\begin{CD}
1 @>>> M @>>> G @>>> B @>>> 1\\
@. @VV\mbox{id}V @VVV @ VV\lambda V @.\\
1 @>>> M @>>> A \wr P @>\bar\alpha >> C\wr P @>>> 1.
\end{CD}
$$
Since $\lambda$ and the identity mapping of $M$ represented by 
the first vertical arrow are monomorphisms, the Short 
Five Lemma implies that the map $G\rightarrow A\wr P$ 
in the middle of the diagram 
is a monomorphism. 
Thus we may consider $G$ as a subgroup of $W=A\wr P$. 
%We may also describe $G$ as
%$$
%G=\{x\in W\mid x\bar\alpha=b\lambda\mbox{ for some }b\in B\}.
%$$
% The group $G$ is said to be the {\em induced extension} defined by the
% homomorphism $\alpha: A\rightarrow B$. 

Induced extensions $G$ satisfy hypothesis~(*), since $\ker\bar\alpha=M$ 
is a normal subgroup of $G$, and since $M$ is the direct product 
$(\ker\alpha)^I$ such that the factors of $(\ker\alpha)^I$ are permuted 
transitively by $G$. One of the main results of \cite{grosskov} is that
the converse is also true.

\begin{theorem}[{\cite[Theorem~4.1]{grosskov}}]\label{induced}
%\begin{enumerate}
%\item The induced extension satisfies the conditions (*).
%\item 
%Conversely, 
If~(*) holds, then 
the map $\alpha:N/K\rightarrow G/M$ given by
$Kx\mapsto Mx$ is a well defined homomorphism
and $G$ is isomorphic to the induced extension defined by $\alpha$.
%\end{enumerate}
\end{theorem}

Let us assume that  $G$ is the induced extension defined by
a homomorphism $\alpha:A\rightarrow B$, so that (*) holds for $G$. 
One of the main objectives of~\cite{grosskov} is to describe  the subgroups $H$ 
of $G$ that satisfy
\begin{equation}\label{subcrit}
HM=G\quad\mbox{and}\quad
H\cap M\cong\prod_{i\in I}(H\cap M)K_i/K_i.
\end{equation}
Such subgroups are called {\em high subgroups} in~\cite{kovacs:max}.
The second assertion of \eqref{subcrit}
implies that $H\cap M$ is not only a direct product, but its direct
product decomposition is inherited from the given direct product 
decomposition of $M$. 
%This is required in case~(2) of Scott's condition.
Condition~\eqref{subcrit} occurs often in the case of primitive 
permutation groups,
%A group $G$ with a subgroup $H$ satisfying~\eqref{subcrit}
%acts on the right coset space $\Omega=[G:H]$. The first equation
%$HM=G$ is equivalent to saying that $M$ is a transitive subgroup of $G$.
%Further, the second equation of~\eqref{subcrit} can be interpreted as saying
%that $\Omega$ admits a direct power decomposition $\Omega=\Gamma^I$ which
%is preserved by $G$. 
since it holds for the point stabilizers of the 
finite primitive groups in several of the O'Nan-Scott classes,  as was 
observed already by Scott 
%the original 
%version of the O'Nan--Scott Theorem given in 
in~\cite{sco:rep}.

The paper \cite{grosskov} describes the subgroups of $G$ that satisfy~\eqref{subcrit}
in terms of subgroups of $N/K$. The following result
plays a pivotal r\^ole in the description of maximal subgroups
in the subsequent paper~\cite{kovacs:max}.

\begin{theorem}[{\cite[Corollary~4.4]{grosskov}}]\label{subgroupsth}
There is a bijection between the conjugacy classes in $G$ of the subgroups 
$H$ that satisfy condition~\eqref{subcrit}
and the conjugacy classes  in $N/K$ of the subgroups $L/K$ that satisfy
$N/K=(M/K)(L/K).$
Under this bijection, complements correspond to complements.
\end{theorem}

Twisted wreath products have played an important part in permutation group
theory, since primitive groups with a non-abelian regular minimal normal 
subgroups (the missing case in the original O'Nan-Scott Theorem) can be
described as twisted wreath products, hence the name `TW-type'.
Kov\'acs and Gross notice
this connection by observing after \cite[Corollary~3.3]{grosskov} 
that if $\alpha:A\rightarrow B$ is such that
 $A$ splits over $\ker\alpha$, then 
the induced extension $G$ can be written as a twisted wreath product
$(\ker \alpha)\,\mbox{\sf twr}\, P$. 
In this case $P$ is a complement of
$M=(\ker\alpha)^I$ in $G$
and
$M$ is a regular  normal subgroup of $G$ in its permutation 
representation on $\Omega=[G:P]$, a crucial observation 
in the study of TW-type primitive groups.

\section{Maximal subgroups}

Determining the subgroup structure of a given finite group has been 
one of the central problems in group theory. The subgroup lattice of a finite
group has a rich structure and its full description is often beyond reach. 
Thus one starts with the maximal subgroups. 

The problem of describing  maximal subgroups of an arbitrary (finite) 
group can be approached
using a simple reduction argument. If $G$ is a group with a normal subgroup
$M$, then the class of maximal subgroups that contain $M$ is in one-to-one
correspondence with the class of maximal subgroups of the quotient $G/M$. 
Hence one really needs to understand core-free maximal subgroups. Core-free
maximal subgroups of $G$ 
are precisely the point stabilizers in faithful primitive permutation
representations of $G$. Hence conjugacy classes of such maximal subgroups
are in a bijective correspondence with equivalence classes of faithful 
primitive permutation representations of $G$. 
% This shows that there is an
% intimate link between obtaining a good understanding of maximal subgroups
% of finite groups and understanding primitive permutation groups.
Thus the problem of describing maximal subgroups of finite groups
can be attacked using the O'Nan--Scott Theorem.
% 
% It was recognized that describing  maximal subgroups of a finite group would
% have to rely on determining maximal subgroups of simple groups. By the early
% 80's, the classification of finite simple groups was accepted to be 
% complete, and attention turned towards applications of the 
% classification theorem, such as the description of maximal subgroups.
% 
% The problem of maximal subgroups 
This problem was tackled in two papers independently, one by
Aschbacher and Scott~\cite{aschsc}, and the second by Laci~\cite{kovacs:max}.
There is significant overlap between the two papers. Laci only learned
about the work of Aschbacher and Scott after  his draft was finished. 
He explains in a footnote added to the second page of his paper:
\begin{quotation}
{\small  After the draft of this paper was completed, I learned that a forthcoming
paper~\cite{aschsc} by Aschbacher and Scott will address the same issues.
$[\ldots]$ While the conclusions naturally have several common components, 
the approaches and expositions differ so much that detailed reconciliation
$[\ldots]$ will not be attempted here. The difficulties involved strongly
suggest that both versions of the story are worth telling. \quad \cite[p.\,115]{kovacs:max}}
\end{quotation}
The  common underlying strategy behind both papers is the one  explained above.
Laci's paper builds on the theory of induced
extensions developed in his  earlier paper~\cite{grosskov} with Gross, 
which allows him to give a more approachable account.

Assuming that $M$ is a minimal normal subgroup of $G$, we are interested 
in the maximal subgroups that do not contain $M$.
In the case of finite groups, 
$M$ is a direct product $T^k$ of simple groups. 
If $T$ is abelian, then so is $M$, and the problem is reduced to
a problem in group cohomology.
%%
%\begin{quotation}
%and the maximal subgroups of $G$ not containing $M$ are precisely 
%the  complements of $M$ in $G$; the number of conjugacy classes of these
%is 0 or the order of the first cohomology group $H^1(G/M,M)$. \quad \cite[p.\,1%14]{kovacs:max}
%\end{quotation}
%The principal part of 
%Laci's paper treats the case when $M$ is non-abelian 
%and non-simple. %The case when $M$ is simple is treated in Section~5.
Hence in the main
part of the paper, $M=T^k$ such that $T$ is a non-abelian simple group, and the $k$ simple 
factors of $M$ form a $G$-conjugacy class. In particular condition (*) 
holds and so, by Theorem~\ref{induced}, $G$ is an
induced extension defined by a natural homomorphism $\alpha:N/K\rightarrow
G/M$ where $K$ is a maximal normal subgroup of $M$ and $N$ is 
its normalizer in $G$.
In particular, there is a correspondence between a certain class
of subgroups of $G$ and the class of subgroups of $N/K$ 
(see Theorem~\ref{subgroupsth}). Laci introduced two types of
subgroups in $G$.
% 
% In a group $G$ that satisfies (*), Laci introduced two types of
% subgroups: {\em high subgroups} and {\em full subgroups}. 
The {\em high subgroups} $H$ defined in \eqref{subcrit} 
may be alternatively defined by requiring that $HM=G$, 
and that either $H\cap M=1$, or (*) holds with $H$, $H\cap M$, $H\cap K$, $H\cap N$ in place
of $G$, $M$, $K$, and $N$. The second class of subgroups, the {\em full subgroups}, are those subgroups $F$
such that $FK=G$ and (*) holds with $F$, $F\cap M$, 
$F\cap K$, $\norm F{F\cap K}$ in  place of $G$, $M$, $K$, $N$.

% High subgroups are the ones that 
% satisfy condition~\eqref{subcrit} which
% appeared in the earlier paper~\cite{grosskov}. 
%To explain the significance of full subgroups, 
%note that
%we require $FK=G$, a stronger condition than $FM=G$. Since $FK=G$, 
%we obtain that 
%$(F\cap M)K=FK\cap M=G\cap M=M$, and hence $(F\cap M)K/K=M/K$. 
%The quotient $(F\cap M)K/K=M/K$ can be viewed as the projection of 
%$F\cap M$ into the direct factor $M/K$ of $M$. 
The condition $FK=G$, in the definition 
of a full subgroup $F$, implies that $(F\cap M)K=M$, and hence that $F\cap M$ projects onto 
$M/K$; and then $FM=G$ yields that it projects onto
each direct factor $M/K_i$. In other words, if 
$F$ is a full subgroup, then $F\cap M$ is a subdirect subgroup
with respect to the direct decomposition $M\cong \prod_i M/K_i$.

The main theorem of \cite{kovacs:max} considers groups $G$ that satisfy (*), with
$M/K$ non-abelian simple and the set $\{K_i\}$ finite, and 
gives a characterization of the maximal subgroups
of $G$ that do not contain $M$.
By~\cite[Lemma~4.2]{kovacs:max}, 
each such  maximal subgroup  of $G$ is either full or high.
Laci split these maximal subgroups into three subfamilies:
\begin{enumerate}
\item[(A)] full maximal subgroups;
\item[(B)] high maximal subgroups $H$ that are not complements of $M$ 
(that is, $H\cap M\neq 1$);
\item[(C)] high maximal subgroups $H$ that are complements of $M$ 
(that is, $H\cap M=1$).
\end{enumerate}
We shall see in a moment why it is necessary to distinguish between 
cases (B) and (C).

In each of these families the number of conjugacy classes of 
maximal subgroups is expressed using  a smaller group. In families (B) and (C), 
this reduction  essentially follows from    
Theorem~\ref{subgroupsth} which states a one-to-one correspondence
between the set of conjugacy classes of high subgroups of $G$ and 
the set of conjugacy classes of supplements of $M/K$ in $N/K$. 
In particular, maximal subgroups in family (B) correspond to 
\begin{enumerate}
\item[(B1)] maximal subgroups $L/K$ of $N/K$ that neither avoid nor contain
$M/K$.
\end{enumerate} 
Describing maximal subgroups in family (C) is a bit more complicated. 
Let $L/K$ be a {\em maximal complement} to $M/K$ in $N/K$ (that is, a maximal subgroup of $N/K$ which complements $M/K$).
The problem is that the corresponding subgroup in $G$, given by 
 Theorem~\ref{subgroupsth},  may not be maximal, as it may 
be contained  in a full subgroup of $G$. One of the main contributions
of the paper~\cite{kovacs:max}  is a sufficient and necessary condition on 
$L/K$ that describes precisely when this situation occurs.
As $N/K=M/K\rtimes L/K$,  the group $N/M$ acts, via the isomorphisms 
$L/K\cong (N/K)/(M/K)\cong N/M$, on $M/K$. Laci proved \cite[Theorem 4.3.c]{kovacs:max} that 
the complement $H$ corresponding in $G$ to $L/K$ is maximal if and only if this action cannot be extended to 
a subgroup of $G/M$ properly containing $N/M$.
Hence maximal subgroups of type (C) correspond to

\begin{enumerate}
\item[(C1)] maximal complements of $M/K$ in $N/K$ such that the corresponding
homomorphism $N/M\rightarrow \aut{(M/K)}$ cannot be extended to a subgroup
of $G/M$ properly containing $N/M$.
\end{enumerate}

The reduction for counting the full maximal 
subgroups is developed in \cite[Section~3]{kovacs:max} with the main result 
being Theorem~3.03. The number of conjugacy classes of full maximal subgroups
is determined in terms of the following set: 
\begin{enumerate}
\item[(A1)] the collection 
of all homomorphisms $\varphi:D\rightarrow \out{M/K}$
where $D$ is a subgroup of $G/M$ minimally containing $N/M$ and the restriction
of $\varphi$ to $N/M$ is equal to the coupling determined by the short exact
sequence 
$$
1\rightarrow M/K\rightarrow N/K\rightarrow N/M\rightarrow 1.
$$
\end{enumerate}

This leads to the main result. 

\begin{theorem}[{\cite[Theorem 4.3]{kovacs:max}}]\label{mainmax}
Suppose that $G$ satisfies (*) with $M/K$  simple and the set $\{K_i\}$  finite.
There is a bijection between the conjugacy classes of maximal subgroups
of $G$ in class (A) and the set of homomorphisms in (A1). Further, there
are bijections between the conjugacy classes of maximal subgroups of $G$ in 
types (B) and (C), and the conjugacy classes of maximal subgroups of $N/M$
of type (B1) and (C1), respectively.
\end{theorem}

The families (A), (B), and (C) of maximal subgroups can 
be interpreted in terms of the O'Nan--Scott classes of the corresponding permutation
representations. Suppose that $G$ is finite with non-abelian, non-simple, minimal
normal subgroup $M$.
Let $H$ be a core-free maximal subgroup of  $G$
and 
 consider the $G$-action on the right cosets of $H$. If $H$ is in
family (A),  then $G$ is a primitive group of SD or CD
 type depending on whether or not $H\cap M$ is simple.
In these cases we also say that the maximal subgroup $H$ is 
of SD or CD type. If $H$ is in family (B), then $G$ is primitive of PA type, 
while if $H$ is a maximal subgroup in family (C), then
$M$ is a regular normal subgroup of $G$ and $G$ can be embedded into the
holomorph of $M$. Here $G$ is either of TW type (when 
$M$ is the unique minimal normal subgroup of $G$) or 
HC type (when $G$ has a second minimal normal subgroup, distinct from $M$).

Full maximal subgroups of SD type  in finite groups 
are further  
studied in the subsequent paper~\cite{kov:sd}. 
The paper uses assumption (*) such that $G$ is finite
and $M$  is a non-abelian, non-simple minimal normal subgroup of $G$.
A maximal subgroup of $G$ with SD type corresponds to 
a homomorphism $G/M\rightarrow \out{M/K}$ in (A1). The main results 
of~\cite{kov:sd} are based on elegant counting arguments in the spirit 
of the quote at the end of Section~\ref{sect1} and describe 
full maximal subgroups of SD type in finite groups.

As observed already in \cite[after Corollary~3.3]{grosskov}, 
finite primitive groups 
with a unique regular non-abelian minimal normal subgroup have a  twisted wreath product structure,
and they are most frequently studied using this structure. 
%A detailed treatment was developed in~\cite{baddeleytw}.
Suppose that $G$ is  a permutation group,  $M$ is a non-abelian 
regular minimal
normal subgroup, and $H$ is a point stabilizer. If $T$ 
is a simple direct factor of $M$, then $M\cong T^k$ and $G$ can be written as
$G=\twr T\varphi H$ where $\varphi$ is the conjugation action
of $Q=\norm HT$ on $T$. By Theorem~\ref{subgroupsth}, the subgroup $H$ 
corresponds
to a complement $L/K$ of $M/K$ in $N/K$.
In the language of twisted wreath products, the subgroups $N/K$, $M/K$, and 
$L/K$ can be identified with $T\rtimes Q$, $T$ and $Q$,
respectively. Now 
the Kov\'acs condition states that 
$H$ is maximal, or equivalently, $G$ is primitive, if and only if 
$Q$ is maximal in $T\rtimes Q$ and the conjugation action of $Q$ cannot
be extended to a larger subgroup of $H$. As $T$ is a finite 
non-abelian simple group,  the condition that $Q$ is maximal in $T\rtimes Q$, 
is equivalent to the condition that $Q$ does not normalize any non-trivial
proper subgroup of $T$, which, in turn, is equivalent to the condition
that conjugation by $Q$ induces a group of automorphisms that contains $\inn T$.
The condition concerning  $\inn T$ is also obtained 
by Aschbacher and Scott in~\cite[Theorem 1(C)(1)]{aschsc}, and conversely 
the `non-extension' property in (C1) can be derived from it. So the papers 
\cite{aschsc,kovacs:max} contain equivalent primitivity conditions for 
twisted wreath type groups.

%\begin{theorem}[Corollary, p.~123 of~\cite{kov:sd}]
%Under the conditions in the previous paragraph, 
%$G$ has at most one conjugacy class of core-free maximal 
%subgroups of SD type, except perhaps when $N/K$ is nearly simple, 
%$G/M$ is soluble, and $m$ is a prime-power divisor of $|\out{M/K}|$. 
%\end{theorem}

%A further significant result of~\cite{kov:sd} is the determination of the numbe%r of isomorphism types of 
%groups of simple diagonal type with a given minimal normal subgroup. 

%\begin{theorem}[Theorem 1 of \cite{kov:sd}]
%Let $T$ be a non-abelian finite simple group and let 
%$k\geq 2$. There is a bijection
%between the permutational isomorphism classes of those primitive groups $G$ 
%of simple diagonal type whose unique minimal normal subgroups are abstractly
%isomorphic to $T^k$ and the conjugacy classes of those subgroups of $\sy k\time%s\out T$ whose projections to $S_k$ are primitive.
%\end{theorem}

Despite not having a published paper devoted to 
primitive groups of TW type, 
Laci's work has, perhaps, left its greatest impact on the theory of
such groups. %In this case, there is a unique minimal normal subgroup
%$M$ which is non-abelian and regular (and hence non-simple), while
%a point stabilizer $H$ is a core-free maximal subgroup complementing $M$.
Laci championed the treatment of these groups as twisted wreath products, and
today it would feel unnatural to treat them in any other way. He and 
Peter F\"orster wrote an
Australian National University Research Report~\cite{kov:tw} on TW groups and,
as witnessed by the bibliographies of~\cite{kovfors,kov:blowup}, he
was working on more. His paper with F\"orster studies conditions
under which the top group $H$ in a twisted wreath product is maximal and
they develop in~\cite[1.1 Theorem]{kovfors} conditions similar to the ones in
(C2) above. It is an important and interesting consequence of their
result~\cite[1.2 Corollary]{kovfors} that if $H$ is such a maximal complement of $M$, then 
$H$ must have a unique minimal normal subgroup $N$ that is non-abelian, and 
a simple factor of $M$ has to occur as a  section in a simple factor
of $N$. This implies the well-known result, that $H$ cannot have 
non-trivial solvable normal subgroups. The same research report contains 
a treatment of the inclusion problem for TW-type groups.

Laci's work on TW-type groups, albeit  formally unpublished, 
had a huge influence on Baddeley's seminal work~\cite{baddeleytw} on 
this topic, as he acknowledged:

\begin{quotation}
{\small  It should be pointed out that a recent research report by 
F\"orster and Kov\'acs~\cite{kov:tw}
contains considerable overlap with our work in \S3 and \S5, and indeed there may
be even more in common between the material in this paper and their
unpublished work.\quad\cite[p.\,548]{baddeleytw}}
\end{quotation}
Section~6 of Baddeley's paper treats the problem of permutational 
isomorphism between abstractly isomorphic twisted wreath products. 
He made Laci's contribution  clear in  a footnote:
\begin{quotation}
{\small  The ideas and results of this section are almost entirely due to 
L.~G.~Kov\'acs. \quad\cite[p.\,568]{baddeleytw}}
\end{quotation}

Laci's work on maximal subgroups  was also highly influential on the 
development of algorithms to determine the conjugacy classes of
maximal subgroups of a finite group.
In fact, a large part of 
the algorithm given in~\cite{eickhulpke} relies on Theorem~\ref{mainmax}.
Later Cannon and Holt~\cite{cannonholt} 
presented essentially an algorithmic version of~\cite{kovacs:max} and in particular
of Theorem~\ref{mainmax}.

\section{Wreath products, blowups, and Wielandt's conjecture}

A lot of work on primitive permutation groups concentrated on understanding
links between primitive groups, non-abelian simple groups
and irreducible representations of finite groups. This information
together with the wreath product construction lead to the solution of many
problems in algebra, number theory and combinatorics. 

The 
paper~\cite{kov:decomp} presents methods for identifying if a primitive
group can be written as a wreath product in product action and provides a means for
counting the possible wreath decompositions. 
Laci introduced, 
as his main tool
for tackling these counting problems,  the concept of a
%Laci introduced a concept, namely the concept of a 
{\em system of product imprimitivity}. A system of product
imprimitivity can be used to detect  embeddings of primitive groups
into wreath 
products in product action in much the same way that
a system of imprimitivity (or block system) can be used to detect embeddings
of transitive groups into wreath products in imprimitive actions.

Primitive wreath products usually have many primitive
subgroups that are not themselves wreath products and, 
for some
applications, detailed information is needed on precisely which subgroups
of a wreath product $H\wr \sy k$ with $H$ primitive on $\Gamma$, are
themselves primitive on $\Gamma^k$. In his seminal paper  \cite{kov:blowup}, Laci
introduced the concept of a `blow-up' of a primitive group and provided
criteria for identifying such subgroups for almost all types of primitive
groups $H$.  Moreover, this led, in 1990, to a classification
\cite{prae:incl} of embeddings of finite
primitive groups into wreath products in product action.

Suppose that $H$ is a primitive group acting on $\Gamma$ with socle $M$. Then 
$W=H\wr \sy k$ can be considered as a permutation group on $\Gamma^k$ in 
its product action. Let $\pi:W\rightarrow \sy k$ be the natural projection
onto $\sy k$, and consider $\pi$ as a permutation representation of $W$.
Let $W_0$ be the stabilizer of $1$ under the representation $\pi$. Then
$W_0$ can be written as a direct product $H\times (H\rtimes \sy {k-1})$. Let
$\pi_0:W_0\rightarrow H$ denote the projection onto the first direct factor.

In the language of~\cite{kov:blowup}, a subgroup $B$ of $W$ is  {\em large} if $B\pi$ is transitive, 
and $(B\cap W_0)\pi_0=H$. Input to the {\em blow-up construction} consists of
a primitive group $G$ on $\Gamma$ with socle $M$ and a large subgroup $B$
of $(G/M)\wr \sy k$ 
with $k\geq 1$. The output $G\uparrow B$, called the blow-up of $G$ by $B$, 
is a permutation group on $\Gamma^k$, 
namely the full inverse image of $B$ under the natural homomorphism
$G\wr\sy k\rightarrow (G/M)\wr\sy k$. 

\begin{theorem}[{\cite[Theorems~1 and~2]{kov:blowup}}]
All blow-ups of a finite primitive permutation group $G$ are primitive 
if and only if the socle of $G$ is not regular. 
If a primitive group $G$ with non-regular socle  is a blow-up, 
then it is a blow-up $G_0\uparrow B$ of a unique $G_0$ which is not itself
a blow-up.
\end{theorem}

To identify which primitive groups are blow-ups, Laci introduced
a  corresponding decomposition concept. %, namely the blow-up decomposition. A 
A {\em blow-up decomposition} of a primitive group $G$ 
with non-regular socle $M$
and point stabilizer $H$ is a direct decomposition of $M$ such that
\begin{enumerate}
\item the direct factors form a $G$-conjugacy class;
\item $H\cap M$ is the direct product of its  intersections with these direct factors.
\end{enumerate}
If $G$ is a primitive group with non-regular socle acting on $\Omega$, then a 
blow-up decomposition of $G$ with $k$ factors 
leads to a bijection $\Omega\rightarrow \Gamma^k$ and 
permutational isomorphism 
$G\rightarrow
 G_0\uparrow B$ where $G_0$ is a primitive group acting on $\Gamma$. 
In other words, such a primitive group that admits
a blow-up decomposition is itself a blow-up of a smaller group.
The primitive groups that admit blow-up decompositions are, in modern 
terminology, 
the groups of PA, CD, HC types and they are blow-ups of 
groups of AS, SD, and HS types, respectively.

Laci's approach that led to the blow-up concept  
also appears in the version of the O'Nan--Scott Theorem presented
in Cameron's  book~\cite{cam}. 
Cameron
divides the finite primitive groups into two large families: the groups in the
first family are called {\em basic groups}, 
and the other family is formed by the {\em non-basic primitive groups}. 
In Cameron's terminology, a primitive group with a non-regular socle is  
basic if and only if it is not a blow-up of a smaller group, while
blow-ups are non-basic.

The philosophy of the blow-up construction is used in the 
paper~\cite{kovfors} jointly written with Peter F\"orster. 
This paper presents an application to Wielandt's problem
on primitive permutation groups. The problem can be stated as follows.
\begin{quotation}
{\small  Suppose that $\pi_1$ and $\pi_2$ are permutation representations of a finite group
$G$ with the same character. Given that $G\pi_1$ is primitive, does it
follow that $G\pi_2$ is primitive? }
\end{quotation}
By obvious reduction, we may assume that $\pi_1$ and $\pi_2$ are 
faithful. Suppose that $G$, $\pi_1$, $\pi_2$ is a counterexample
to Wielandt's conjecture and let $A$ and $B$ be point stabilizers for
the  representations $\pi_1$ and $\pi_2$, respectively. Then
$A$ is maximal while $B$ is not, and we may assume that $A$ and $B$ are core-free.
If $\soc G$ is not simple, the argument of the paper shows that $G$, 
considered as a transitive group with point stabilizer $A$, is a primitive
group of PA-type. As discussed above, this implies that 
$G$ is a blow-up of an almost simple group $G_0$ and it is shown in 
the paper that the almost simple group $G_0$ is also a counterexample
to the conjecture. In fact the paper shows how to construct the 
set of all faithful counterexamples, given the set of faithful almost simple
counterexamples. Hence the conjecture is reduced to the case of 
almost simple groups. 

The Wielandt conjecture was proven to be false by Guralnick
and Saxl~\cite{guralnicksaxl} who presented  the first
almost simple counterexamples. 
Later more counterexample emerged in the work of 
Breuer~\cite{breuer}.
% Later more counterexample surged in the works of 
% Breuer~\cite{breuer}.

\section{Other work on permutation groups}

%\subsection*{Minimal faithful degree}
The {\em minimal faithful degree} $\mu(G)$ of a finite group $G$ is the 
the size of the smallest set on which $G$ can be represented faithfully 
as a permutation group. Certainly, by Cayley's theorem, $\mu(G)\leq |G|$, but it can be considerably smaller.
It is interesting to know (for example, in computational applications) about
the relation, if any, between $\mu(G)$ and $\mu(G/N)$ where $N$ is a normal subgroup of $G$. 
In general, $\mu(G/N)$ can be much larger than $\mu(G)$. In fact Neumann~\cite{neu86} 
gave examples for which $\mu(G/N)>c^{\mu(G)}$ with $c=2^{0.25}$, while Holt and Walton~\cite{holwal}
showed that $\mu(G)$, for arbitrary finite groups $G$, satisfies $\mu(G/N)\leq c^{\mu(G)-1}$ with $c=4.5$.
On the other hand, if $G/N$ has no nontrivial abelian normal subgroup (for example, if $N$ is the soluble radical of $G$), then Laci and the first 
author~\cite{kovpra00} 
showed that $\mu(G/N)\leq\mu(G)$. They proved the same inequality if $G/N$ is elementary abelian  \cite[p.\,284]{kovpra89}, and 
conjectured that $\mu(G/N)\leq\mu(G)$ should hold whenever $G/N$ is abelian. 

Although this conjecture remains unresolved, Theorem 2 of~\cite{kovpra00} 
shows the following (where the last assertion is by Franchi~\cite{fra11}).

\begin{theorem} 
For a potential counterexample $(G,N)$ with both $\mu(G)$ and $|G|$ 
minimal, $G$ must be
a directly indecomposable $p$-group for some prime $p$, $N$ must be 
the derived subgroup $G'$,
$\mu(G/N) = \mu(G) + p$, and $G$ has no abelian maximal
subgroup.
\end{theorem}

The problem  above, for $G/N$ abelian, motivated the study~\cite{kovpra89} of permutation groups with nontrivial 
abelian quotients. Suppose that $G$ is a permutation group on $n$ points and that, for some prime divisor $p$ of $|G|$, 
a Sylow $p$-subgroup moves exactly $kp$ points. Then (see \cite[Theorem]{kovpra89}) the largest abelian $p$-quotient of 
$G$ has order at most $p^k$, with equality if and only if $G$ is the direct product of its largest $p'$-constituent and
its transitive non-$p'$-constituents, and the possibilities for the latter are explicitly listed. 
As a corollary: if $G$ is not perfect, that is, if $G\ne G'$, then for 
some prime divisor $p$ of $|G/G'|$, the order $|G/G'|\leq p^{n/p}$, and the groups $G$ for which equality holds are
described. The cyclic group of order $30$ represented as a permutation group on $10$ points, demonstrates that 
this upper bound may hold for more than one but not all prime divisors of $|G/G'|$. The result about $\mu(G/N)$ for $G/N$ 
elementary abelian is then deduced. Another immediate consequence  is that
$|G/G'|\leq 3^{n/3}$ for all permutation groups $G$ on $n$ points, and so $|G|\leq 3^{n/3}$ if $G$ is an abelian permutation group on $n$ points.
This upper bound is markedly different from the situation for primitive groups $G$ where it was proved in \cite{aschgur} that 
$|G/G'|\leq n$ and that $n$ must be prime if equality holds.

There are several other papers by Laci on permutation groups to be noted. In the first \cite{cknp} from 1985, the authors show that
a transitive permutation group $P$ of order a power of a prime $p$ has at least $(p|P|-1)/(p+1)$ fixed-point-free elements,
and analyses the possible structures of groups attaining the bound. A second paper~\cite{kovnew88}, written in collaboration with M.~F.~Newman,
concerns the number of generators for finite nilpotent transitive groups. 
They prove that there is a constant $c$ such that each 
nilpotent transitive permutation group on $n$ points can be generated by a set of $cn(\log n)^{-1/2}$ elements,
and on the other hand they show that, for each prime $p$ there is a constant $c_p$ such that, for each $p$-power $n$,
there is a transitive $p$-group on $n$ points which cannot be generated by $c_pn(\log n)^{-1/2}$ elements. 
The upper bound is now known to hold for all transitive permutation 
groups~\cite{bkr,lmm}.

Laci, with G.~R.~Robinson~\cite{kr} obtained an exponential 
upper bound $5^{n-1}$ on the number of conjugacy classes of a 
permutation group of degree $n$. 
Despite many improvements over the years, 
the best currently known general upper bound is still of 
this form with $5^{1/3}$ in place of $5$~\cite{maroti}.

\def\cprime{$'$}

\end{document}